\input amstex
\documentstyle{amsppt}
\magnification=\magstep1 \baselineskip=12pt \hsize=6truein
\vsize=9truein

\topmatter

\title
On positive solutions to semi-linear conformally invariant equations
on locally conformally flat manifolds
\endtitle

\author Jie Qing and David Raske
%\footnote{Department of Mathematics UC Santa Cruz, Santa Cruz, CA
%95064, qing{\@}ucsc.edu, raske{\@}ucsc.edu}
\endauthor

\leftheadtext{Compactness and Existence} \rightheadtext{Jie Qing and
David Raske}

\address Jie Qing, Department of Mathematics, UC Santa Cruz, Santa Cruz, CA
95064 \endaddress \email qing{\@}ucsc.edu\endemail
\address David Raske, Department of Mathematics, UC Santa Cruz, Santa Cruz, CA
95064 \endaddress \email gadfly{\@}ucsc.edu\endemail

\abstract In this paper we study the existence and compactness of
positive solutions to a family of conformally invariant equations on
closed locally conformally flat manifolds. The family of conformally
covariant operators $P_\alpha$ were introduced via the scattering
theory for Poincar\'{e} metrics associated with a conformal manifold
$(M^n, [g])$. We prove that, on a closed and locally conformally
flat manifold with Poincar\'{e} exponent less than $\frac
{n-\alpha}2$ for some $\alpha \in [2, n)$, the set of positive
smooth solutions to the equation
$$
P_\alpha u = u^\frac {n+\alpha}{n-\alpha}
$$
is compact in the $C^\infty$ topology. Therefore the existence of
positive solutions follows from the existence of Yamabe metrics and
a degree theory.
\endabstract

\endtopmatter

\document

\head 1. Introduction \endhead

In a recent paper of Graham and Zworski [GZ] a meromorphic family of
conformally covariant operators associated with a Riemannian
manifold $(M^n, g)$ was introduced via scattering theory for a
Poincar\'{e} metric associated with $(M^n, [g])$. For a metric $g\in
[g]$ on $M^n$, the scattering operator $S(z)$ is a meromorphic
family of peusdo-differential operators on $M$ for $\text{Re}(z)>
n/2$ of order $2\text{Re}(z) - n$ with the principal symbol
$$
\sigma(S(z)) = 2^{n-2z}\frac {\Gamma(n/2-z)}{\Gamma(z-n/2)}
\sigma((-\Delta_{g})^{z-n/2}).
$$
The scattering operator is conformally covariant in the sense that
$$
S(z)[e^{2\Upsilon}g] = e^{- z\Upsilon}S(z)[g]e^{(n-z)\Upsilon}.
$$
At the special values $z=(n/2)+k$,
$$
\text{Res}_{z = \frac n2 + k} S(z) = c_k P_{2k},
$$
where $c_k = (-1)^{k+1}[2^{2k-1}k!(k-1)!]^{-1}$ and $P_{2k}$ is the
conformally invariant powers of the Laplacian introduced in [GJMS].
For instance,
$$
P_2[g] = -\Delta_g + \frac {n-2}{4(n-1)}R[g] \tag 1.1
$$
is the conformal Laplacian and
$$
P_4 [g] = (-\Delta)^2 - \text{div}_g((\frac {(n-2)^2 +
4}{2(n-1)(n-2)} R[g] - \frac 4{n-2}Ric[g])d) + \frac {n-4}2 Q[g]
\tag 1.2
$$
is the Paneitz operator [P], where $R[g]$ is the scalar curvature,
$Ric[g]$ is the Ricci curvature, and
$$
Q[g] = = - \frac 1{2(n-1)} \Delta R + \frac {n^3 - 4n^2 + 16n -
16}{8(n-1)^2(n-2)^2}R^2 - \frac {4}{(n-2)^2}|Ric|^2 \tag 1.3
$$
is the so-called $Q$-curvature.

The conformal Laplacian governs the transformation of the scalar
curvature under the conformal change of metrics
$$
P_2[g] u = \frac {n-2}{4(n-1)}R[u^\frac 4{n-2} g] u^\frac
{n+2}{n-2}; \tag 1.4
$$
while the Paneitz operator $P_4$ governs the transformation of $Q$
curvature
$$
P_4[g] u = \frac {n-4}2 Q[u^\frac 4{n-4}g]u^\frac {n+4}{n-4}. \tag
1.5
$$
The well-known Yamabe problem in differential geometry is to find a
metric of constant scalar curvature in a given class of conformal
metrics, that is, to solve the Yamabe equation
$$
P_2[g] u = u^{\frac {n+2}{n-2}} \tag 1.6
$$
for some positive function $u$ on a given manifold $(M^n, g)$
($n\geq 3$). The affirmative resolution to the Yamabe problem was
given in [Sc1] after other notable works [Ya] [Tr] [Au]. Recent
developments in the study of conformal geometry and conformally
invariant partial differential equations have created increasingly
interests in the higher order generalizations of Yamabe problem as
are described in [CY1] (references therein). One may ask whether
there is a metric of constant $Q$ curvature in a given conformal
class of metrics, that is, to solve the Paneitz-Branson equation
$$
P_4[g] u = u^\frac {n+4}{n-4} \tag 1.7
$$
for some positive smooth function $u$ on a given manifold $(M^n, g)$
($n\geq 5$) [P] [Br] [DHL] [DMA].

In our previous paper [QR], we considered the Paneitz-Branson
equation (1.7) and were able to obtain the compactness of metrics of
constant $Q$ curvature in a given conformal class on a locally
conformally flat manifold. But it has been a very challenging
problem to find a positive solution to the higher order
Paneitz-Branson equation. We discovered it is useful to consider a
smooth family of equations
$$
P_\alpha[g] u =  u^\frac {n+\alpha}{n-\alpha} \tag 1.8
$$
where
$$
P_\alpha [g] = 2^\alpha \frac {\Gamma(\frac \alpha 2)}{\Gamma(-\frac
\alpha 2)} S(\frac {n+\alpha}2)[g] \tag 1.9
$$
and $\alpha \in [2, n)$.

On a locally conformally flat manifold $(M^n, g)$ with positive
Yamabe constant, there is the family of conformally covariant
operators $P_\alpha$ associated with the hyperbolic metric whose
conformal infinity is $(M^n, [g])$. Readers are referred to a paper
of Patterson and Perry [PP] for more detailed discussions on the
scattering theory for the conformally compact hyperbolic manifolds.
If the Poincar\'{e} exponent of the holonomy representation of the
fundamental group $\pi_1(M)$ is less than $\frac {n-\alpha}2$, then
there is a conformally covariant integral operator $I_\alpha$, which
is a right inverse to $P_\alpha$. The existence of such integral
operators allow us to consider the conformally invariant integral
equation
$$
u = I_\alpha [g] (u^\frac {n+\alpha}{n-\alpha}). \tag 1.10
$$
We observed that the new moving plane method introduced by Chen, Li
and Ou in [CLO] (please also see [CY2]) can be adopted to study the
conformally invariant integral equation (1.10) and allows us to
prove the following

\proclaim{Theorem 1.1} Suppose that $(M^n, g)$ is a locally
conformally flat manifold with positive Yamabe constant and
Poincar\'{e} exponent less than $\frac {n-\alpha}2$ for some $\alpha
\in [2, n)$. And suppose that $(M^n, g)$ is not conformally
equivalent to the standard round sphere. Then there exists a
constant $C=C(n, \alpha, k)$ such that, for any smooth positive
solution $u$ to (1.10), we have
$$
\|u\|_{C^k(M)} + \|\frac 1u\|_{C^k(M)} \leq C. \tag 1.11
$$
\endproclaim

Consequently, based on the degree theory described in [Sc2], we can
prove the existence of positive solutions to both equations (1.8)
and (1.10).

\proclaim{Theorem 1.2} Suppose that $(M^n, g)$ is a locally
conformally flat manifold with positive Yamabe constant and
Poincar\'{e} exponent less than $\frac {n-\alpha_0}2$ for some
$\alpha_0 \in [2, n)$, then, for any $\alpha \in [2, \alpha_0]$,
there exists a positive smooth function $u$, which solves both
equations (1.8) and (1.10).
\endproclaim

In particular,

\proclaim{Theorem 1.3} Suppose that $(M^n, g)$ ($n\geq 5$) is a
locally conformally flat manifold with positive Yamabe constant and
Poincar\'{e} exponent less than $\frac {n-4}2$, then there exists a
positive smooth solution $u$ to the Paneitz-Branson equation (1.7).
\endproclaim

The organization of this paper is as follows. We will first describe
in Section 2 the conformally covariant integral operator $I_\alpha$.
Then we will discuss its relation to the conformally covariant
pseudo-differential operator $P_\alpha$ in Section 3. In Section 4
we adopt the new moving plane method to derive an apriori estimate
in the form of a convexity theorem for solutions to (1.10). In
Section 5 we will introduce the rescaling method for the conformally
invariant integral equations (1.10), obtain the compactness and
prove the existence.

\head  2. Conformally covariant operators on LCF manifolds
\endhead

We start with a discussion of the scattering operator and an
introduction of the integral operators on locally conformally flat
manifolds. We first recall from Proposition 4.1 in [PP] that on the
hyperbolic space $H^{n+1}$ in the half space model the scattering
operator is
$$
S^0(z) = 2^{n -2z} \frac {\Gamma(\frac n2 -z)}{\Gamma(z - \frac n2)}
(-\Delta)^{2z - n}. \tag 2.1
$$
Hence
$$
P^0_\alpha = 2^\alpha \frac {\Gamma(\frac \alpha 2)}{\Gamma(-\frac
\alpha 2)}S^0(\frac {\alpha+n}2) = (-\Delta)^{\frac \alpha 2}. \tag
2.2
$$

Suppose that $(M^n, g)$ is a closed locally conformally flat
manifold with positive Yamabe constant. According to Schoen-Yau
[SY], the developing map from the universal cover of $M^n$ into the
round sphere $(S^n, g_1)$ is injective. The deck transformation
group of this covering becomes a Kleinian group $\Gamma$. The image
of the developing map is the set $\Omega(\Gamma)$ of ordinary points
for the Kleinian group $\Gamma$, and $(M^n, [g])$ is conformally
equivalent to $\Omega(\Gamma)/\Gamma$. Therefore we may consider the
universal covering
$$
\phi: (\Omega(\Gamma), \tilde g) \to (M, g) \tag 2.3
$$
and write $\tilde g = \phi^* g = \tilde \eta^2 g_1$. Using
stereographic projection
$$
\psi: \hat\Omega \subset R^n \longrightarrow \Omega (\Gamma)\subset
S^n \tag 2.4
$$
with respect to some point in $\Omega(\Gamma)\subset S^n$, we may
write $\hat g  =  \psi^* \tilde g= \hat \eta^2 g_0$, where $g_0$ is
the Euclidean metric and $\hat \eta = (\tilde\eta \circ \psi)(\frac
2{1+|x|^2})$. Thus we unfolded a locally conformally flat manifold
$(M^n, g)$ into a conformally flat manifold $(\hat \Omega, \hat g)$.

One way to understand the scattering operator for the conformally
compact hyperbolic manifold $H^{n+1}/\Gamma$ is to lift everything
to the hyperbolic space $H^{n+1}$. Upon using the geodesic defining
function corresponding to the choice of the metric $g\in [g]$ on
$M$, we can equivalently consider the lifting from $M$ to $\hat
\Omega\subset R^n$. Let $u$ be a function on $M$ and $\hat u = u
\circ \phi\circ \psi$. Then
$$
P_\alpha[g] u = \hat P_\alpha [\hat g]\hat u|_F, \tag 2.5
$$
where $F$ is a fundamental domain for $\Gamma$. By the conformally
covariant property, we have
$$
\hat P_\alpha[\hat g] \hat u = \hat \eta^{- \frac {n+\alpha}2} \hat
P_\alpha [g_0] \hat \eta^{\frac {n-\alpha}2} \hat u, \tag 2.6
$$
where clearly $\hat P_\alpha [g_0] = P^0_\alpha = (-\Delta)^\frac
\alpha 2$ due to the uniqueness of solutions to the Poisson
equations. Notice that $\hat \eta^{\frac {n-\alpha}2} \hat u \in
L^1(R^n)$ (cf. Corollary 2.2 in [QR]).

We now turn to consider an integral operator on the Euclidean space
$R^n$
$$
\hat I_\alpha (\hat u)(x) = \int_{R^n} \frac {c(n, \alpha)}{|x -
y|^{n-\alpha}} \hat u(y)dy, \quad \alpha \in (0, n),
$$
for a function $\hat u$ on $R^n$, where
$$
c(n, \alpha) = \frac {\Gamma(\frac {n-\alpha}2)}{\pi^{\frac n2
-\alpha}\Gamma(\frac \alpha 2)}. \tag 2.6
$$
It is well known that
$$
(-\Delta)^\frac \alpha 2 \hat I_\alpha (\hat u) = \hat u, \quad
\alpha \in (0, n). \tag 2.7
$$

Then, by the conformally covariant property, we have the integral
operator, on $(\Omega(\Gamma), \tilde g)$,
$$
\tilde I_\alpha[\tilde g](\tilde u)(\tilde x) = \int_{S^n} \tilde
K(\tilde x, \tilde y)\tilde u(\tilde y) dV_{\tilde g}(\tilde y),
\tag 2.8
$$
where
$$
\tilde K(\tilde x, \tilde y) = \tilde \eta^{- \frac {n-\alpha}2}
(\tilde x)(\frac {1+|x|^2}2)^\frac {n-\alpha}2 \frac {c(n,
\alpha)}{|x-y|^{n-\alpha}}(\frac {1+|y|^2}2)^\frac {n-\alpha}2\tilde
\eta^{-\frac {n-\alpha}2} (\tilde y). \tag 2.9
$$
Notice that $S^n \setminus \Omega (\Gamma)$ is of Hausdorff
dimension less than $\frac {n-2}2$ for a locally conformally flat
manifold $(M^n, g)$ with positive Yamabe constant. Hence we have
$$
\tilde I_\alpha[\tilde g](\tilde u)(\tilde x) =
\int_{\Omega(\Gamma)} \tilde K(\tilde x, \tilde y)\tilde u(\tilde y)
dV_{\tilde g}(\tilde y). \tag 2.10
$$
To push down the operator to $(M^n, g)$ we simply notice that we can
always lift a function on $M$ to a periodic function on its covering
$\Omega(\Gamma)$, which is,
$$
\tilde u (\tilde x) = u(\phi(\tilde x)).
$$
So we define, for a fundamental domain $F\subset \Omega(\Gamma)$ and
$\tilde x\in F$,
$$
\aligned I_\alpha (u)(p) & = \int_{\Omega(\Gamma)} \tilde K(\tilde
x, \tilde y) \tilde u(\tilde y)dV_{\tilde g}(\tilde y) =
\sum_{\gamma\in\Gamma} \int_{\gamma F}\tilde K(\tilde x, \tilde y)
u(q)dV_{\tilde g}(\tilde y) \\ & = \sum_{\gamma\in\Gamma}
\int_{F}\tilde K(\tilde x, \gamma \tilde y) u(q) dV_{\tilde
g}(\tilde y) = \int_M K(p, q) u(q) dV_g(q),\endaligned \tag 2.11
$$
where $\phi(\tilde x) = p$, $\phi(\gamma \tilde y) = q$ and
$$
K(p, q)  = \sum_{\gamma\in\Gamma} \tilde K(\tilde x, \gamma \tilde
y). \tag 2.12
$$

\proclaim{Theorem 2.1} Suppose that $(M^n, [g])$ is a locally
conformally flat manifold with positive Yamabe constant and that the
Poincar\'{e} exponent of the holonomy representation is less than
$\frac {n-\alpha}2$. Then the integral operator $I_\alpha$ given in
the above is a well-defined conformally covariant operator of
bi-degree $(\frac {n-\alpha}2, \frac {n+\alpha}2)$ on the space
$C^\infty(M)$. Moreover
$$
P_\alpha[g] I_\alpha [g] u = u, \tag 2.13
$$
for any function $u$ on $M$.
\endproclaim

\demo{Proof} Following the above discussion we need to show that the
kernel $K(p, q)$ is well-defined and that $I_\alpha$ is independent
of the choice of a fundamental domain $F$. From (2.9) and
$$
\eta (\gamma \tilde x) = |\gamma'(\tilde x)|^{-1}\eta (\tilde x).
\tag 2.14
$$
we then have
$$
\aligned & \sum_{\gamma\in \Gamma} \tilde K(\tilde x, \gamma\tilde
y)\\
= \tilde \eta^{- \frac {n-\alpha}2}(\tilde x)(\frac
{1+|x|^2}2)^\frac {n-\alpha}2  & \tilde \eta^{-\frac {n-\alpha}2}
(\tilde y)\sum_{\gamma\in \Gamma}\frac {c(n, \alpha)}{|x - \gamma
y|^{n-\alpha}}(\frac {1+|\gamma y |^2}2)^\frac
{n-\alpha}2|\gamma'|^\frac {n-\alpha}2 (\tilde y).
\endaligned
$$
Notice that there is some constant $K >0$ such that
$$
\frac {c(n, \alpha)}{|x - \gamma y|^{n-\alpha}}(\frac {1+|\gamma y
|^2}2)^\frac {n-\alpha}2 \leq K, \forall \ \gamma \in
\Gamma\setminus \{1\}. \tag 2.15
$$
Therefore, the convergence of $K(p, q)$ depends only on the
convergence of the Poincar\'{e} series
$$
\sum_{\gamma\in \Gamma} |\gamma'|^\frac {n-\alpha}2. \tag 2.16
$$
We know the above Poincar\'{e} series is convergent since the
Poincar\'{e} exponent is less than $\frac {n-\alpha}2$. To prove
that the kernel function $K(x, y)$ does not depends on $F$, we only
need to show that
$$
\tilde K(\gamma\tilde x, \gamma\tilde y) = \tilde K(\tilde x, \tilde
y). \tag 2.17
$$
We recall (1.3.2) from [Ni] that
$$
|\gamma x - \gamma y| =|\gamma'(\tilde x)|_e^\frac 12 |x
-y||\gamma'(\tilde y)|_e^\frac 12, \tag 2.18
$$
where $|\gamma'|_e$ is the norm under the Euclidean metric. Then the
above (2.17) follows from (2.18) and the following
$$
|\gamma'(\tilde x)| = \frac {1+|x|^2}{1+|\gamma x|^2}|\gamma'(x)|_e.
\tag 2.19
$$
One may find a similar calculation in the proof of Proposition 1.1
in [CQY].

To prove (2.13) we simply calculate, for a function $u\in
C^\infty(M)$ and $v = I_\alpha [g] (u)$, in the light of (2.5),
$$
P_\alpha [g] v = \hat \eta^{-\frac {n+\alpha}2}(-\Delta)^\frac
\alpha 2 \hat \eta^{\frac {n-\alpha}2} \hat v
$$
where, by (2.11) and (2.9),
$$
\hat \eta^{\frac {n-\alpha}2} \hat v  = \int_{R^n} \frac {c(n,
\alpha)}{|x-y|^{n-\alpha}} (\frac 2{1+|y|^2})^\frac{n + \alpha}2
\tilde\eta^\frac {n+\alpha}2(\psi( y)) \hat u(y) dy.
$$
Thus the proof is completed.
\enddemo

\head 3. Conformally invariant equations
\endhead

Suppose that $(M^n, g)$ is a locally conformally flat manifold with
positive Yamabe constant, and that the Poincar\'{e} exponent of the
holonomy representation of the fundamental group $\pi_1(M)$ is less
than $\frac {n-\alpha}2$ for some $\alpha \in [2, n)$. We consider
the conformally invariant integral equation
$$
u (p) = I_\alpha (u^\frac {n+\alpha}{n-\alpha}) = \int_M K(p, q)
u^\frac {n+\alpha}{n-\alpha} (q) dV_g(q) \tag 3.1
$$
for a positive function $u(p)\in C^\infty (M)$, where $K(p, q)$ is
define as in Theorem 2.1 in the previous section. Let
$$
\hat v (x) = u(\phi(\psi(x)))\tilde\eta^\frac{n-\alpha}2(\psi(x))
(\frac 2{1+|x|^2})^\frac {n-\alpha}2 \in C^\infty(\hat \Omega). \tag
3.2
$$
Then we have
$$
\hat v (x) = \int_{\hat \Omega} \frac {c(n, \alpha)}{|x
-y|^{n-\alpha}}\hat v^\frac {n+\alpha}{n-\alpha} (y) dy. \tag 3.3
$$

Equation (2.13) shows that a positive solution $u$ to the integral
equation (3.1) solves the pseudo-differential equation
$$
P_\alpha u = u^\frac {n+\alpha}{n-\alpha} \quad \text{on $M$} \tag
3.4
$$
and a positive solution to (3.4) solves the integral equation (3.1)
when $P_\alpha$ has a trivial kernel. Thus we have

\proclaim{Theorem 3.1} Suppose that $(M^n, [g])$ is a closed locally
conformally flat manifold with positive Yamabe constant, and that
the Poincar\'{e} exponent of the holonomy representation of the
fundamental group $\pi_1(M)$ is less than $\frac {n-\alpha}2$ for
some $\alpha \in [2, n)$. Then a positive smooth function $u$ solves
the pseudo-differential equation
$$
P_\alpha u = u^\frac{n+\alpha}{n-\alpha} \quad \text{on $M$},
$$
if and only if it solves the integral equation
$$
u = I_\alpha (u^\frac{n+\alpha}{n-\alpha})  \quad \text{on $M$}.
$$
\endproclaim

\demo{Proof} We need to prove that, when the Hausdorff dimension of
the limit set of the holonomy representation of the fundamental
group $\pi_1(M)$ is less than $\frac {n-\alpha}2$, the operator
$P_\alpha$ has a trivial kernel. (2.13) says $P_\alpha$ has a right
inverse $I_\alpha$ when the Poincar\'{e} exponent of the holonomy
representation of the fundamental group $\pi_1(M)$ is less than
$\frac {n-\alpha}2$. Hence $P_\alpha$ is surjective and its
co-kernel is trivial. Therefore $P_\alpha$ has a trivial kernel
since $P_\alpha$ for $\alpha\in [2, n)$ is self-adjoint from [GZ]
(see also an easier proof in [FG]).
\enddemo

For the scalar curvature and conformal Laplacian $P_2$, the Yamabe
constant is defined as
$$
Y(M, [g]) = \inf\{ \frac{\int_M \phi P_2[g]\phi dV_g}{(\int_M
\phi^{\frac {2n}{n-2}}dV_g)^\frac {n-2}n}: \phi \in C^\infty(M), \
\phi > 0 \}, \tag 3.5
$$
which is a conformal invariant of $(M, [g])$. We may similarly
define the Yamabe constant of order $\alpha$ as follows:
$$
Y_\alpha(M, [g]) = \inf\{ \frac{\int_M \phi P_\alpha[g]\phi
dV_g}{(\int_M \phi^{\frac {2n}{n-\alpha}}dV_g)^\frac {n-\alpha}n}:
\phi \in C^\infty(M), \ \phi > 0\}. \tag 3.6
$$
This is a family of conformal invariants of $(M, [g])$ since
$$
\frac{\int_M \phi P_\alpha[g]\phi dV_g}{(\int_M \phi^{\frac
{2n}{n-\alpha}}dV_g)^\frac {n-\alpha}n} = \frac{\int_M (e^{-\frac
{n-\alpha}2\Upsilon} \phi) P_\alpha[e^{2\Upsilon}g](e^{-\frac
{n-\alpha}2\Upsilon}\phi) dV_{e^{2\Upsilon}g}}{(\int_M (e^{-\frac
{n-\alpha}2\Upsilon}\phi)^{\frac
{2n}{n-\alpha}}dV_{e^{2\Upsilon}g})^\frac {n-\alpha}n}.
$$
We know that the first eigenvalue $\lambda_1(P_2)$ of $P_2$ is
positive if and only if the Yamabe constant is positive. We observed

\proclaim{Theorem 3.2} Suppose that $(M^n, g)$ is a locally
conformally flat manifold with positive Yamabe constant. Suppose
that the Poincar\'{e} exponent of the holonomy representation of the
fundamental group $\pi_1(M)$ is less than $\frac {n-\alpha}2$ for
some $\alpha \in [2, n)$. Then the Yamabe constant of order $\alpha$
is positive.
\endproclaim

\demo{Proof} We first claim that the first eigenvalue
$\lambda_1(P_\beta)$ for any $\beta \in [2, \alpha]$ of $P_\beta$ is
positive. This is because we know $\lambda_1(P_2)$ is positive.
Then, due to the continuity of the operator $P_\beta$ with respect
to $\beta$, $\lambda_1(P_\beta)$ has to be zero before getting to
negative. But, from the proof of the previous Theorem 3.1, we know
$P_\beta$ is injective for all $\beta\in [2, \alpha]$. Hence
$\lambda_1(P_\beta)$ remains positive in $[2, \alpha]$ from $\beta =
2$. Now, we simply apply the Sobolev inequality
$$
\| u \|_{L^\frac{2n}{n-\alpha}(M)} \leq C \|u\|_{W^{\frac \alpha 2,
2}(M)}
$$
to conclude $Y_\alpha > 0$ since $\int_M u P_\alpha u dV$ is
equivalent to $\|u\|_{W^{\frac \alpha 2, 2}(M)}$ by the definition
of $P_\alpha$ in this case.
\enddemo

\head 4. Convexity \endhead

To obtain a priori estimates for positive solutions to the integral
equation (3.1). We use the same approach as in [Sc2] [QR]. We will
first establish a convexity theorem. The new moving plane method
introduced by Chen, Li and Ou in [CLO], which was motivated by a
work in [CY2], is what it takes for us to study the solutions to the
class of integral equations (3.1). Let us state and prove the
following convexity theorem.

\proclaim{Theorem 4.1} Suppose that $(M^n, [g])$ is a locally
conformally flat manifold with positive Yamabe constant, and that
the Hausdorff dimension of the limit set of the holonomy
representation of the fundamental group $\pi_1(M)$ is less than
$\frac {n-\alpha}2$ for some $\alpha \in [2, n)$. In addition we
assume that $(M^n, [g])$ is not conformally equivalent to the
standard round sphere. Let $u$ be a positive smooth solution to
(3.1). Then every round ball in $(M, [g])$ is geodesically convex
with respect to the metric $u^\frac 4{n-\alpha} g$.
\endproclaim

\demo{Proof} Fix any point $\tilde x_0\in \partial B\subset
\Omega(\Gamma)$ and consider the stereographic projection with
respect to the antipode of $\tilde x_0$ on $\partial B$. Hence
$$
\psi^{-1}(\partial B) = \{x\in R^n: x_n = 0\}. \tag 4.1
$$
It suffices to prove that, with respect to each $\tilde x_0\in
\partial B$, the hyperplane $\{ x_n=0\}$ is geodesically convex in
the metric $\hat u^\frac 4{n-\alpha} g_0$, where
$$
\hat v (x) = u \circ\phi\circ\psi (x) \tilde\eta^\frac {n-\alpha}2
(\psi (x)) (\frac 2{1+|x|^2})^\frac {n-\alpha}2 \tag 4.2
$$
is a positive solution to
$$
\hat v(x) = \int_{\hat \Omega}\frac {c(n, \alpha)}{|x -
y|^{n-\alpha}} \hat v^\frac {n+\alpha}{n-\alpha}(y)dy. \tag 4.3
$$
%That is to prove that
%$$
%\frac{\partial \hat v}{\partial x_n}|_{x_n=0} < 0. \tag 4.4
%$$
%Notice that
%$$
%\hat L = \psi^{-1}(L(\Gamma)) \subset \{x\in R^n: x_n < 0\} \tag 4.5
%$$
%is a compact subset in the lower half space of $R^n$.
Now let us start the moving plane method introduced in [CLO]. Let
$$
\Sigma_\lambda = \{x\in R^n: x_n > \lambda\} \quad \text{and} \
S_\lambda = \{x\in R^n: x_n = \lambda\}.
$$
We consider the reflection with respect to the hyperplane
$S_\lambda$
$$
x_\lambda = (x_1, x_2, \cdots, x_{n-1}, 2\lambda - x_n), \forall \
x\in R^n \tag 4.4
$$
and define
$$
\hat v_\lambda (x) = \hat v (x_\lambda). \tag 4.5
$$
Then, following a simple calculation, we have
$$
\aligned \hat v(x) & - \hat v_\lambda(x) \\= c(n, \alpha)
\int_{\Sigma_\lambda}(|x - y|^{-n+ \alpha} & - |x -
y_\lambda|^{-n+\alpha}) (\hat v^\frac{n+\alpha}{n-\alpha} (y) - \hat
v_\lambda^\frac{n+\alpha}{n-\alpha} (y)) dy.
\endaligned \tag 4.6
$$
%Notice that
%$$
%|x - y| < |x - y_\lambda|, \forall \ x\in \Sigma_\lambda. \tag 4.13
%$$
Hence it follows from Lemma 2.1 in [CLO] that, for some constant
$C$,
$$
\hat u(x) - \hat u_\lambda(x) < C \int_{\Sigma_\lambda^-} |x -
y|^{-n+\alpha} \hat u^\frac {2\alpha}{n-\alpha}(y) (\hat u (y) -
\hat u_\lambda (y))dy, \tag 4.7
$$
where
$$
\Sigma_\lambda^-= \{x \in \Sigma_\lambda: \hat v_\lambda (x) < \hat
v(x)\}. \tag 4.8
$$
We recall the classic Hardy-Littlewood-Sobolev inequality
$$
|\int_{R^n}\int_{R^n} f(x)|x-y|^{\alpha - n} g(y)dydx| \leq c(n,
\alpha, s)\|f\|_{L^r(R^n)} \|g\|_{L^s(R^n)}, \tag 4.9
$$
where
$$
\frac 1r + \frac 1s = 1 + \frac \alpha n,
$$
$f\in L^r(R^n)$ and $g\in L^s(R^n)$. Therefore, if let
$$
\aligned  f = (\hat v - \hat v_\lambda)^{\mu -
1}\chi_{\Sigma_\lambda^-}, \quad & g = \hat v^\frac
{2\alpha}{n-\alpha}(\hat v - \hat v_\lambda
)\chi_{\Sigma_\lambda^-},
\\ \nu = \frac n\alpha, \quad \mu > \frac n{n-\alpha}, \quad & \frac
1s = \frac 1\nu + \frac 1\mu, \quad \frac 1r + \frac 1\mu =
1,\endaligned \tag 4.10
$$
from (4.7) we have
$$
\|\hat v - \hat v_\lambda\|^\mu_{L^\mu(\Sigma_\lambda^-)} \leq c(n,
\alpha, \mu) C \|\hat v - \hat
v_\lambda\|_{L^\mu(\Sigma_\lambda^-)}^\mu \|\hat v\|^{\frac \alpha
n}_{L^{\frac{2n}{n-\alpha}}(\Sigma_\lambda^-)}, \tag 4.11
$$
which implies
$$
c(n, \alpha, \mu) C \|\hat v\|^{\frac \alpha
n}_{L^{\frac{2n}{n-\alpha}}(\Sigma_\lambda^-)} \geq 1, \tag 4.12
$$
if
$$
\hat v - \hat v_\lambda\in L^\mu(\Sigma_\lambda^-). \tag 4.13
$$
We observe that
$$
\hat v \in L^{\frac {n+\alpha}{n-\alpha}}(R^n), \tag 4.14
$$
which can be proven in a way similar to Corollary 2.2 in [QR].
(4.14) implies (4.13) since $\hat L$ is a compact subset of $R^n$.

We notice that, when $\lambda$ is very large, $\Sigma_\lambda$
corresponds to a very small round ball $B_\epsilon \subset B\subset
\Omega(\Gamma) \subset S^n$. Therefore, when $\lambda$ is very
large,
$$
|\Sigma_\lambda^-| = 0,
$$
which implies $\Sigma_\lambda^-$ is empty when $\lambda$ is very
large. Thus we have
$$
\hat v (x) \leq \hat v_\lambda (x), \ \forall \ x\in \Sigma_\lambda
\tag 4.15
$$
when $\lambda$ is sufficiently large. This gets the moving plane
started. Next we show that we can move the hyperplane down as long
as it does not touch the singular set $\hat L$, which is the image
of the limit set under the stereographic projection $\psi$. Suppose
that
$$
\lambda_0 = \inf\{\lambda: \hat v(x) \leq v_\lambda (x), \forall \
x\in \Sigma_\lambda\}
$$
and
$$
\Sigma_{\lambda_0}\bigcap \hat L  = \emptyset.
$$
Notice that $\hat L$ is not empty because $(M^n, [g])$ is not
conformally equivalent to the round sphere.  Then, by (4.15), we
know
$$
\hat v(x) < \hat v_{\lambda_0} (x), \forall \ x\in
\Sigma_{\lambda_0}. \tag 4.16
$$
Hence, as observed in [CLO],
$$
\lim_{\lambda\to\lambda_0} |\Sigma_\lambda^-| = 0. \tag 4.19
$$
Therefore there is some small number $\epsilon >0$, such that
$$
\Sigma_{\lambda_0 -\epsilon} \bigcap \hat L = \emptyset \tag 4.20
$$
and
$$
c(n, \alpha, \mu) C \|\hat v\|^{\frac \alpha
n}_{L^{\frac{2n}{n-\alpha}}(\Sigma_\lambda^-)} < \frac 12, \forall \
\lambda_0 - \epsilon \leq \lambda \leq \lambda. \tag 4.21
$$
Thus $\Sigma_\lambda^-$ has to be empty and
$$
\hat v(x) \leq \hat v_\lambda (x), \forall \ x\in \Sigma_\lambda
$$
for $\lambda \geq \lambda_0 - \epsilon$, which contradicts with the
definition of $\lambda_0$. Particularly we have, for $\lambda = 0$,
$$
\hat v (x) < \hat v_0(x), \forall x\in \Sigma_0 \tag 4.22
$$
and
$$
\frac{\partial \hat v}{\partial x_n}|_{x_n = 0} = (n-\alpha)c(n,
\alpha) \int_{\Sigma_0} \frac {y_n(\hat u^\frac {n+\alpha}{n-\alpha}
(y) - \hat u_0^\frac {n+\alpha}{n-\alpha}(y))}{|x
-y|^{n-\alpha+2}}dy < 0. \tag 4.23
$$
So the convexity is proven.
\enddemo

\head 5. A priori estimate and existence \endhead

In this section we derive a priori estimates for solution to the
integral equation (3.1) and pseudo-differential equation (3.4) on a
locally conformally flat manifold $(M^n, [g])$, when the
Poincar\'{e} exponent of the holonomy representation of its
fundamental group $\pi_1(M)$ is less than $\frac {n-\alpha}2$. The
idea will be the same as in [Sc2] and [QR], which is to use the
convexity to eliminate possible blow-ups. Due to the global nature
of the equations we in the following describe, for instance, the
rescaling method for integral equation (3.1). For that, we first
establish the $C^{k,\theta}$ estimates based on the $L^\infty$
bounds.

\proclaim{Lemma 5.1} Suppose that $(M^n, [g])$ is a locally
conformally flat manifold $(M^n, [g])$ such that the Poincar\'{e}
exponent of the holonomy representation of its fundamental group
$\pi_1(M)$ is less than $\frac {n-\alpha}2$ for some $\alpha \in [2,
n)$. And suppose that $u$ is a positive smooth solution to the
equation (3.1) on $M$. Then, for $\theta\in (0, 1)$,
$$
\|u \|_{C^{k, \theta} (M)} \leq C (\|u\|_{L^\infty(M)}). \tag 5.1
$$
\endproclaim

\demo{Proof} We may unfold the function $u$ and consider $\hat v$ on
$\hat \Omega$, which satisfies
$$
\hat v  = \int_{R^n} \frac {c(n, \alpha)}{|x-y|^{n-\alpha}} \hat
v^\frac {n+\alpha}{n-\alpha}(y) dy. \tag 5.2
$$
Then the estimate follows from standard elliptic theory.
\enddemo

Next suppose that $u$ is a positive smooth solution to (3.1). Let
$p_0 \in M$ be a fixed point in $M$ and let
$$
\zeta (x) = \phi(\psi( \tilde\eta^\frac 2{n-\alpha} (0) x)):
B_{2\delta}(0) \to M, \ \text{and} \ \zeta_\lambda (x) = \zeta
(\frac x{\lambda^\frac 2{n-\alpha}}). \tag 5.3
$$
Then if let
$$
v_\lambda (x) = \frac 1\lambda u(\zeta_\lambda(x)),
$$
we have, for $x\in B_{\lambda^{\frac 2{n-\alpha}}\delta}$,
$$
\aligned v_\lambda(x) = \int_M \frac 1\lambda K(\zeta_\lambda (x),
q)  & u^\frac {n+\alpha}{n-\alpha} dV_g\\ = \int_{B_{\lambda^\frac
2{n-\alpha}\delta}} \frac 1{\lambda^2} K(\zeta_\lambda (x),
\zeta_\lambda (y)) v_\lambda^\frac {n+\alpha}{n-\alpha}(y)
dV_{g_\lambda}(y) & + \frac 1\lambda \int_{M \setminus \zeta
(B_\delta)}K(\zeta_\lambda (x), q) u^\frac {n+\alpha}{n-\alpha}
dV_g,
\endaligned \tag 5.4
$$
where $g_\lambda = \lambda^{-\frac 4{n-\alpha}}g$. Meanwhile
$$
\hat v_\lambda (x) = \int_{R^n}\frac {c(n,
\alpha)}{|x-y|^{n-\alpha}} \hat v_\lambda^\frac
{n+\alpha}{n-\alpha}(y) dy. \tag 5.5
$$
Therefore, from (2.9), we have

\proclaim{Lemma 5.2} In the above, for any fixed $\Lambda
> 0$,
$$
\left\{\aligned \frac 1{\lambda^2} K(\zeta (\frac x{\lambda^\frac
2{n-\alpha}}), \zeta (\frac y{\lambda^\frac 2{n-\alpha}})) &
\rightrightarrows \frac {c(n, \alpha)}{|x-y|^{n-\alpha}} \\
dV_{g_\lambda}(y) & \rightrightarrows |dy|^2, \endaligned\right.
\quad \text{on $B_\Lambda(0)$} \tag 5.5
$$
as $\lambda \to \infty$.
\endproclaim

On the other hand,  there is a priori integral bound for any
positive smooth solutions as follows:

\proclaim{Lemma 5.3} Suppose $u$ is a positive smooth solution to
the equation (3.4). Then, for a constant $C=C(M, g)$,
$$
\int_M u^{\frac {n+\alpha}{n-\alpha}}dv_g < C(M, g). \tag 5.6
$$
\endproclaim

\demo{Proof} One simply integrates the equation (3.4) and gets
$$
\int_M u^\frac {n+\alpha}{n-\alpha} dv_g = \int_M P_\alpha u dv_g =
\int_M u P_\alpha 1 dv_g \tag 5.7
$$
due to the fact that $P_\alpha$ is self-adjoint (cf. [GZ] [FG]).
Hence
$$
\int_M u^\frac {n+\alpha}{n-\alpha} dv_g \leq \max_M |P_\alpha[g]1|
(\int_M u^\frac {n+\alpha}{n-\alpha} dv_g)^\frac
{n-\alpha}{n+\alpha} \text{vol}_g(M)^\frac {2\alpha}{n+\alpha}. \tag
5.8
$$
Therefore
$$
\int_M u^\frac {n+\alpha}{n-\alpha} dv_g \leq \max_M
|P_\alpha[g]1|^\frac {n+\alpha}{2\alpha}\text{vol}_g(M). \tag 5.9
$$
\enddemo

Now, let us consider a sequence of smooth solutions $\{u_k\}$ to
(3.1) such that
$$
\lambda_k = \max_M u_k  = u_k (p_k) \to \infty.
$$
Let $v_k$ be the rescaled function as defined in the above. Then we
have, for any $\Lambda >0$, via a priori estimates in Lemma 5.1 and
(5.5),
$$
v_k(x) \rightrightarrows v(x) \quad \text{in $C^1(B_\Lambda)$}, \tag
5.10
$$
where $v \in C^1_{loc}(R^n)$, at least for some subsequence of
$\{u_k\}$. Moreover it follows from the above Lemma 5.3 that the
second term in (5.4) always converges to zero as $\lambda_k\to
\infty$. Hence
$$
v (x) =  \int_{R^n} \frac {c(n, \alpha)}{|x-y|^{n-\alpha}} v^{\frac
{n+\alpha}{n-\alpha}} (y) dy. \tag 5.11
$$

Therefore we have

\proclaim{Theorem 5.4} Suppose that $(M^n, [g])$ is a locally
conformally flat manifold with positive Yamabe constant such that
the Hausdorff dimension of the limit set of the holonomy
representation of its fundamental group is less than $\frac
{n-\alpha}2$ and $(M, [g])$ is not conformally equivalent to the
standard round sphere. Then there exists a constant $C=C(n, \alpha,
k)$ such that, for any positive smooth solution $u$ to the integral
equation (3.1),
$$
\|u\|_{C^k(M)} + \|\frac 1u\|_{C^k(M)}  \leq C(n, \alpha, k). \tag
5.12
$$
\endproclaim

\demo{Proof} We first use the above rescaling method to derive a
priori $L^\infty$-bound. Assume otherwise there are a sequence of
smooth solutions as the above sequence $\{u_k\}$. Then we obtain a
smooth solution $v$ to the integral equation on $R^n$  as in (5.11).
The classification given in [CLO], which generalizes the result in
[Ln] and [WX], shows $v^\frac 4{n-\alpha} |dx|^2$ is isometric to
the standard round sphere (may not be the unit round sphere though).
Hence the sufficiently large ball $B_K = \{x\in R^n: |x| < K\}$ has
a concave boundary. Therefore, in the light of (5.10), such $v$
could not exist because of the convexity theorem in the previous
section (see similar arguments in [Sc2] [QR]). Higher order
estimates follows from standard elliptic theory.

Next we use the fact that the Yamabe constant of order $\alpha$ is
positive by Theorem 3.2 to prove the positive lower bound. Assume
otherwise, there is a sequence of smooth solutions $\{u_k\}$ such
that
$$
\inf_M u_k = u_k(p_k) \to 0. \tag 5.13
$$
Then, because of the $L^\infty$ estimates in the first step, there
is a subsequence $\{u_k\}$ converging strongly, say in $C^2(M)$, to
a solution $u$ with $u(p_0) = 0$ for some $p_0\in M$, which implies
that $u\equiv 0$. But,
$$
\frac {\int_M u_k P_\alpha u_k dV_g}{(\int_M u_k^\frac
{2n}{n-\alpha} dV_g)^\frac {n-\alpha}n} = (\int_M u_k^\frac
{2n}{n-\alpha} dV_g)^\frac \alpha n \to 0, \tag 5.14
$$
which contradicts with the assumption that $Y_\alpha(M, [g]) > 0$.
Thus the proof is completed.
\enddemo

Finally we adopt the degree theory approach from [Sc2] to prove the
existence of solutions to (3.1), therefore existence of solutions to
(3.4). We consider, for $\theta\in (0, 1)$,
$$
\Omega_\Lambda = \{u\in C^{2, \theta}(M): \|u\|_{C^{2, \theta}(M)} +
\|\frac 1u\|_{C^{2, \theta}(M)} < \Lambda, u > 0 \}
$$
and the map
$$
F_\alpha = u - I_\alpha(u^\frac {n+\alpha}{n-\alpha}):
\Omega_\Lambda \to C^{2, \theta}(M).
$$
From the elliptic theory we know that $F_\alpha = Id +
\text{compact}$ and we may define the Leray-Schauder degree (cf.
[N]) of $F_\alpha$ in the region $\Omega_\Lambda$ with respect to
$0\in C^{2, \theta}(M)$, denoted by $\text{deg}(F_\alpha,
\Omega_\Lambda, 0)$, provided that $0\notin F_\alpha(\partial
\Omega_\Lambda)$.

\proclaim{Theorem 5.5} Suppose that $(M^n, g)$ is a locally
conformally flat manifold with positive Yamabe constant. Suppose
that the Hausdorff dimension of the limit set of the holonomy
representation of the fundamental group $\pi_1(M)$ is less than
$\frac {n-\alpha_0}2$ for some $\alpha_0\in [2, n)$. Then, for each
$\alpha \in [2, \alpha_0]$, there exists a positive smooth solution
to both the integral equation (3.1) and the pseudo-differential
equation (3.4).
\endproclaim

\demo{Proof} Suppose that $(M^n, g)$ is not conformally equivalent
to the standard round sphere. Otherwise we know a standard round
metric is the solution to both (3.1) and (3.4). Then, by Theorem 5.4
in the above, there is a number $\Lambda
> 0$ such that $0 \notin F_\alpha (\partial \Omega_\Lambda)$ for
any $\alpha\in [2, \alpha_0]$. And the homotopy invariance of the
degree tells us that
$$
\text{deg}(F_\alpha, \Omega_\Lambda, 0) = \text{deg}(F_2,
\Omega_\Lambda, 0).
$$
for all $\alpha \in [2, \alpha_0]$. Hence, due to [Sc2],
$$
\text{deg}(F_\alpha, \Omega_\Lambda, 0) = \text{deg}(F_2,
\Omega_\Lambda, 0) = -1. \tag 5.15
$$
Therefore there is a positive solution $u\in C^{2, \theta}(M)$ to
the integral equation (3.1). Thus, by elliptic regularity theory,
$u$ is smooth and $u$ also solves the pseudo-differential equation
(3.4).
\enddemo

\vskip 0.1in \noindent {\bf References}:

\roster \vskip0.1in
\item"{[Au]}" T. Aubin, The scalar curvature, differential geometry
and relativity, (Cahen and Flato, eds.), Reidel, Dordrecht 1976.

\vskip0.1in
\item"{[Br]}" T. Branson, Group representations arising from
Lorentz conformal geometry, J. Func. Anal. 74 (1987), 199-291.

%\vskip0.1in
%\item"{[CGY1]}" S.-Y. A. Chang, M. Gursky and P. Yang, An equation of
%Monge-Ampere type in conformal geometry and 4-manifolds of
%positive Ricci curvature, Annals of Math. 155 (2002), 709-787.

%\vskip 0.1in
%\item"{[CGY2]}" S.-Y. A. Chang, M. Gursky and P. Yang,
%A conformally invariant Sphere theorem in four dimension, Preprint
%2002.

%\vskip 0.1in
%\item"{[CHY]}" S.-Y. A. Chang, F. Hang and P. Yang, On a class of
%locally conformally flat manifold, Preprint 2003.

%\vskip0.1in
%\item"{[CQY1]}" S.-Y. A. Chang, J. Qing and P. Yang, On the
%compactification of a class of locally conformally flat
%4-manifolds, Invent. Math. 142 (2000), no. 1, 65--93.

\vskip0.1in
\item"{[CQY]}" S.-Y. A. Chang, J. Qing and P. Yang, On the finiteness
of Kleinian groups in general dimensions, J. Reine Angew. Math., 571
(2004), 1-17.

\vskip0.1in
\item"{[CY1]}" S.-Y. A. Chang and P. Yang, Nonlinear differential
equations in conformal geometry, Proceedings of ICM 2002.

\vskip0.1in
\item"{[CY1]}" S.-Y. A. Chang and P. Yang, On uniqueness of solution
of a n-th order differential equaton in conformal geometry, Math.
Res. Letter, 4 (1997) 91-102.

\vskip0.1in
\item"{[CLO]}" W. Chen, C. Li and B. Ou, Classification of
solutions for an integral equation, to appear in Comm. Pure and
Appl. Math.

\vskip0.1in
\item"{[DHL]}" Z. Djadli, E. Hebey and M. Ledoux, Paneitz type
operators and applications, Duke Math. J. 104(2000) no. 1, 129-169.

\vskip0.1in
\item"{[DMA]}" Z. Djadli, A. Malchiodi and M. Ahmedou, Prescribing
a fourth order conformal invariant on the standard sphere, Part II -
blow-up analysis and applications, Preprint 2001.

%\vskip0.1in
%\item"{[GNN]}" B. Gidas, W.-M. Ni and L. Nirenberg, Symmetry and
%related properties via maximum principle, Comm. Math. Phys. 68
%(1979), no. 3, 209-243.

%\vskip0.1in
%\item"{[G]}" Maria Gonzalez, Singular sets of a class of locally
%conformally flat manifolds, Preprint 2004.

\vskip 0.1in
\item"{[GJMS]}" C. R. Graham, R. Jenne, L. Mason and G. Sparling,
Conformally invariant powers of the Laplacian, I. Existenece, J.
London Math. Soc. (2) 46 (1992) 557-565.

\vskip 0.1in
\item "{[GZ]}" C. R. Graham and M. Zworski, Scattering matrix
in conformal geometry, Invent. Math. 152(2003) no. 1, 89-118.

%\vskip0.1in
%\item"{[Gu]}" Matt Gursky, The Weyl functional, de Rham
%cohomology, and K\"{a}hler-Einstein metrics, Annals. Math. 148
%(1998), 315-337.

\vskip0.1in
\item"{[Ln]}" C.-S. Lin, A classification of solutions of a conformally
invariant fourth order equation in $R^n$, Comment. Math. Helv. 73
(1998), 206-231.

\vskip 0.1in
\item"{[Ni]}" P. Nicholls; The Ergodic theory of discrete groups. London Math. Soc.
Lecture Note Ser., 143. Cambridge Univ. Press, Cambridge, 1989.

\vskip0.1in
\item"{[N]}" L. Nirenberg, Topics in nonlinear functional
analysis, Courant Institute Publication 1973-74.

\vskip0.1in
\item"{[P]}" S. Paneitz, A quadratic conformally covariant
differential operator for arbitrary pseudo-Riemannian manifolds,
Preprint 1983.

\vskip0.1in
\item"{[PP]}"  S. Patterson and P. Perry, The divisor of Selberg's
zeta function for Kleinian groups, Duke Math. J. 106 (2001), no. 2,
pp 321-390.

%\vskip0.1in
%\item"{[Q]}" Jie Qing, A Priori Estimates for Positive
%Solutions of Semi-linear Elliptic Systems. J. Partial Differential
%Equations, Vol.1, No.2, Series A(1988), pp 61-70.

\vskip0.1in
\item"{[QR]}" J. Qing and D. Raske, Compactness for conformal metrics
with Constant $Q$ curvature on locally conformally flat manifolds,
to appear in Cal. Var. PDE.

\vskip0.1in
\item"{[Sc1]}" R. Schoen, Conformal deformation of a Riemmanian
metric to constant scalar curvature, J. Diff. Geom. 6 (1984),
479-495.

\vskip0.1in
\item"{[Sc2]}" R. Schoen, On the number of constant scalar
curvature metrics in a conformal class, Differential Geometry,
311-32o, Pitman Monogr. Surveys Pure Appl. Math., 52, Longman Sci.
Tech., Harlow, 1991.

\vskip0.1in
\item"{[SY]}" R. Schoen and S.T. Yau, Conformally flat manifolds,
Kleinian groups, and scalar curvature, Invent. Math. 92 (1988),
47-71.

\vskip0.1in
\item"{[Tr]}" N. Trudinger, Remarks concerning the conformal
deformation of Riemannian structures on compact manifolds, Ann.
Scuola Norm. Sup. Pisa Cl. Sci (3) 22 (1968), 265-274.

\vskip0.1in
\item"{[WX]}" J. Wei and X. Xu, Classification of
solutions of higher order conformally invariant equations, Math.
Ann. 313 (1999), 207-228.

\vskip0.1in
\item"{[Ya]}" H. Yamabe, On a deformation of Riemannian structures
on compact manifolds, Osaka J. Math. 12 (1960), 21-37.

\endroster
\enddocument

\enddocument